\documentclass[12pt]{amsproc}%
\usepackage{amsmath}
\usepackage{amsfonts}
\usepackage{amssymb}
\usepackage{graphicx}%
\providecommand{\U}[1]{\protect\rule{.1in}{.1in}}
\newtheorem{theorem}{Theorem}[section]

\newtheorem{proposition}[theorem]{Proposition}
\newtheorem{example}[theorem]{Example}

\def\U{{\mathcal{U}}}

\begin{document}
\author{Peter A. Loeb}
\address{Department of Mathematics, University of Illinois, Uabana, IL 61801. e-mail: PeterA3@AOL.com}
\dedicatory{This article is dedicated to the memory of \ Wim Luxemburg. His friendship and
encouragement remain a source of investigative courage. Great thanks are also
due to Renming Song for many helpful conversations.}\title[Martin Boundary]{An Intuitive Approach to the Martin Boundary}
\date{November 1, 2019}

\begin{abstract}
An intuitive probabilistic alternative for the construction of the Martin
boundary is presented along with a construction of maximal representing
measures for positive harmonic functions.

\end{abstract}
\maketitle

\section{Introduction}

Robert Martin's 1941 boundary construction in \cite{Martin} is now a
fundamental tool in potential theory and probability theory. (See
\cite{Brelot} and \cite{Doob}.) We suggest here an intuitive probabilistic
alternative to Martin's Green's function construction. It begins with the
author's result with M. Insall and M. Marciniak in \cite{Compactifications}
showing that a compactifying boundary is formed by equivalence classes of
points not infinitely close to standard points in the nonstandard extension of
a metric (or even a regular) space. An extension in \cite{Hausdorff
Compactifications} shows that any Hausdorff compactification can be formed in
this way. This raises the following question: What is an equivalence relation
using probability theory for Martin's boundary?\footnote{Presented, in part,
at the May 2019 Calgary meeting of the Statistical Society of Canada.}

The equivalence relation presented here produces a compactification that
coincide with Martin's for important examples. We also show, using
\cite{Loeb-Boundary}, that the resulting construction leads to the correct
representing measures for positive harmonic functions. This approach
\textquotedblleft looking inside\textquotedblright\ a domain only makes sense
when one can speak of points that are neither points of an existing boundary
nor points in a compact subset of the domain.

I am grateful to the editors of this memorial volume for the opportunity to
acknowledge my lasting debt to Wim Luxemburg and to describe, with the hope of
beginning a continuing conversation, the results obtained so far.

\section{General Compactifications}

In this section, we review joint work with Insall and Marciniak in
\cite{Compactifications} and \cite{Hausdorff Compactifications}. We use basic
concepts of Robinson's \cite{Robinson} nonstandard analysis that are discussed
in Appendix C of \cite{Loeb-Real} and more deeply in the beginning of
\cite{LW}. While the results of this section extend to a regular topological
space, we restrict our attention here to a noncompact metric space $(W,d)$. A
compactification of $W$ is a compact space containing $W$ as dense subspace.

Recall that if $r$ is a nonstandard real number a finite distance from $0$,
then the standard part of $r$ is denoted by $\operatorname{st}r$. It is the
unique standard real number infinitely close to $r$. Also, we let {}$^{\ast
}\mathbb{N}_{\infty}$ denote the set of unlimited natural numbers, i.e.,
$^{\ast}\mathbb{N}\backslash\mathbb{N}$.

Fix a nonstandard extension $~^{\ast}W$ of $W$. A point $y\in\,^{\ast}W$ is
called \textbf{near-standard} if there is a standard point point $x\in W$ with
{}$^{\ast}d(y,x)\simeq0$. That is, the distance from $x$ to $y$ is
infinitesimal. If $y\in\,^{\ast}W$ is not near-standard, then we say that $y$
is \textbf{remote} in $^{\ast}W$. For example, the remote points in the
nonstandard extension of the open unit disk are the points with distance from
the origin infinitely close to $1$. By the main result in
\cite{Compactifications} (summarized in Chapter 5 of \cite{LW}), given an
equivalence relation on the remote points of $^{\ast}W$, the equivalence
classes form the boundary points of a compactification of $W$. By
\cite{Hausdorff Compactifications}, every Hausdorff compactification can be
formed this way. The following example extends even to topological spaces that
are just regular.

\begin{example}
[\textbf{Stone-\v{C}ech}]%
\begin{rm}%
Let $\mathcal{F}$ be the set of all bounded continuous real-valued functions
on $W$. The Stone-\v{C}ech compactification is produced by the equivalence
relation that sets remote points $x$ and $y$ equivalent if and only if for all
$f$ in $\mathcal{F}$,
\[
^{\ast}f{}(x)-~^{\ast}f{}(y)\simeq0.
\]%
\end{rm}%

\end{example}

\section{A Probabilistic Equivalence Relation}

Let $W$ be an open connected domain in Euclidean space or manifold suitable
for potential theory and Brownian motion. Let $n\mapsto K_{n}$, $n\in
\mathbb{N}$, be a compact exhaustion of $W$. That is, each $K_{n}$ is a
compact set contained in the interior of $K_{n+1}$, and $\cup_{n\in\mathbb{N}%
}K_{n}=W$. We may assume that each $K_{n}$ is the closure of a connected
region that is regular for the Dirichlet problem. By Robinson's compactness
criterion, a point $x$ is remote in~$^{\ast}W$ if and only if it is
outside$~^{\ast}K_{n}~$\ for every $n\in\mathbb{N}$.

For all $n\in\mathbb{N}$ and $z\in W\backslash K_{n}$, let $\rho_{z}^{n}$ be
the conditional probability measure on $\partial K_{n}$ given by the exit
distribution of a Brownian particle for $W\backslash K_{n}$. That is,
$\rho_{z}^{n}$ is the restriction to $\partial K_{n}$ of harmonic measure for
$W\backslash K_{n}$ but normalized to total mass $1$. The function
$z\mapsto\rho_{z}^{n}$ extends to the points of {}$^{\ast}W\backslash{}^{\ast
}K_{n}$.

Given remote points $x$ and $y$, we set%
\[
f_{n}(x,y):=\operatorname{st}\left(  \left\vert {}^{\ast}\rho_{x}^{n}-{}%
^{\ast}\rho_{y}^{n}\right\vert \left(  {}^{\ast}\partial K_{n}\right)
\right)  .
\]
Note that $f_{n}(x,y)$ is a real-valued sequence. We call remote points $x$
and $y$ equivalent, and write $x\sim y$, if that sequence has limit $0$. That
is,%
\[
x\sim y\Longleftrightarrow\lim_{n\rightarrow\infty}f_{n}(x,y)=0.
\]
There is still a great deal of freedom in choosing an exhaustion $n\mapsto
K_{n}$.

\section{Preserving an Existing Boundary}

A topological boundary of $W$ may be reproduced by the above equivalence
relation. Recall, for example, the celebrated result of Hunt and Wheeden
\cite{Hunt-Wheeden}, showing that the topological boundary of a Lipschitz
domain is the Martin boundary. In any case, each remote point $x$ will be
infinitely close to a unique standard point $\operatorname{st}x$ in the
topological boundary. The following condition is clear.

\begin{proposition}
A topological boundary of $W$ will be produced by the equivalence relation
$\sim$ if and only if for each pair of remote points $x$ and $y$,%
\[
x\sim y\Longleftrightarrow\operatorname{st}x=\operatorname{st}y.
\]

\end{proposition}

\begin{example}
[Unit Disk]\label{UnitDisk}%
\begin{rm}%
For each $n\in\mathbb{N}$, let $K_{n}$ be the closed disk of radius $1-1/n$
centered at the origin. Two remote points $x$ and $y$ each have a standard
part on the unit circle, and $x\sim y$ if and only if $\operatorname{st}%
x=\operatorname{st}y$. Therefore, the corresponding compactification is the
closed unit disk.%
\end{rm}%

\end{example}

\begin{example}
[Cut Unit Disk]\label{CutUnitDisk}%
\begin{rm}%
Let $W$ be the open unit disk from which the interval $(0,1)$ in the $x$-axis
has been removed. Given a positive $\varepsilon\simeq0$, the remote points
$1/2+\varepsilon i$ and $1/2-\varepsilon i$ have the same standard part in the
complex plane, but they are not equivalent. The equivalence relation replaces
the part of the topological boundary formed by interval $(0,1]$ on the
$x$-axis with two similar but distinct intervals.%
\end{rm}%

\end{example}

\section{Two Exotic Examples}

Here are two more examples where the equivalence relation produces the Martin boundary.

\begin{example}
\label{Spheres}%
\begin{rm}%
Let $S$ be the topological boundary of the sphere in $\mathbb{R}^{3}~$of
radius $2$ centered at the point $(0,0,2)$. Let $T$ be the solid sphere of
radius $1$ centered at $(0,0,1)$. Clearly, $T$ is contained inside $S$, and
the boundary of $T$ intersects $S$ just at the origin. Let $W$ be the open
region between $T$ and $S$. As compact sets $K_{n}$ fill $W$, they must
surround much of $T$ above the origin. We may assume that the bottom boundary
of sets $K_{n}$ for large $n$ are annular $2$-dimensional regions. It follows
that the boundary produced by the equivalence relation replaces the origin in
the topological boundary of $W$ with a ring of points. This coincides with the
Martin boundary for $W$.%
\end{rm}%

\end{example}

\begin{example}
\label{Comb}%
\begin{rm}%
Start with the open square in $\mathbb{R}^{2}$ given by $0<x<1$, $0<y<1$.
Assume that for each $n\in\mathbb{N}$, the following interval has been removed
from the square.%
\[
x=1/n\text{,~~~~~}0<y\leq1-1/n\text{.}%
\]
Let $W$ be the resulting region. Suppose $(\xi,\eta)$ is a remote point in
{}$^{\ast}W$ with $\xi$ strictly between $1/\omega$ and $1/(\omega+1)$, where
$\omega$ is in {}$^{\ast}\mathbb{N}_{\infty}$. Then the path of any Brownian
particle starting at the point that exits $W$ from the boundary of a standard
compact set $K_{n}$ must contain points for which the $x$-coordinate is
infinitesimal and the $y$ coordinate is infinitely close to $1$. It follows
that the boundary produced by the equivalence relation intersects the $y$-axis
at the point $(0,1)$. This coincides with the Martin boundary for $W$.%
\end{rm}%

\end{example}

\section{Representing Measures}

Fix $x_{0}\in W$. Let $\mathcal{H}^{1}$ be space of all positive harmonic
functions $h$ on $W$ with $h(x_{0})=1$. The set $\mathcal{H}^{1}$ is convex
and compact with respect to the topology of uniform convergence on compact
subsets of $W$, i.e., the \textbf{ucc topology}. Each $h\in\mathcal{H}^{1}%
\;$is represented by\ a unique probability measure on the extreme elements of
$\mathcal{H}^{1}$. We will use the following construction of that measure from
\cite{Loeb-Boundary}. It is an early application, discussed at the 1974
Oberwolfach Conferences on Potential Theory, of the general measure
construction later published in \cite{Loeb}. A standard, weak-limit
construction of representing measures, established using the nonstandard
construction, can be found in \cite{StandardConstruction}.

Fix $\gamma\in~^{\ast}\mathbb{N}_{\infty}$. Let $U$ denote the internal
interior of $K_{\gamma}$; let $C$ denote $\partial K_{\gamma}$. Recall that
$C$ is an internal Dirichlet regular\ boundary of $U$. For \ each $x\in U$,
let $\mu_{x}$ denote the internal harmonic measure for $x$ on $C$. That is,
for each $x\in U$ and each internally continuous $g$ on $C$, the map%
\[
x\mapsto\int_{C}g(s)d\mu_{x}(s)
\]
gives the value at $x$ of the internal harmonic extension of $g$ from $C$ to
$U$.

We next fix a hyperfinite partition of $C$ consisting of internal Borel
subsets. We assume the partition is so fine that for every $h\in
\mathcal{H}^{1}$,{}$^{\ast}h$ has infinitesimal variation on each set in the
partition. We denote by $\{A_{i}\}$ an internal indexed subfamily of the
partition such that for each index $i$, $\mu_{x_{0}}\left(  A_{i}\right)  >0$,
and $\mu_{x_{0}}\left(  C\backslash\cup_{i}A_{i}\right)  =0$, whence for each
$x\in U$, $\mu_{x}\left(  C\backslash\cup_{i}A_{i}\right)  =0$.

For each index $i$, fix $y_{i}\in A_{i}$. The~{}$^{\ast}\mathbb{R}$-valued map
on $U$ given by
\[
z\mapsto\dfrac{\mu_{z}(A_{i})}{\mu_{x_{0}}(A_{i})}%
\]
is an internal harmonic function in $U$. It is the solution for the function
that is $1$ on $A_{i}$ and $0$ on the rest of $C$, then normalized to be $1$
at $x_{0}$. The real-valued function $\operatorname{uccst}\left(  \mu_{\cdot
}(A_{i})/\mu_{x_{0}}(A_{i})\right)  $ given by
\[
x\mapsto\operatorname{st}\left(  \dfrac{\mu_{x}(A_{i})}{\mu_{x_{0}}(A_{i}%
)}\right)  ,~~~x\in W
\]
is its standard part in $\mathcal{H}^{1}$ with respect to the the ucc topology.

Given any standard $h\in\mathcal{H}^{1}$ and $x\in W$,%

\begin{align*}
h(x) &  =\int_{C}{}^{\ast}\!h(y)\,d\mu_{x}(y)\\
&  \simeq\sum_{i}{}^{\ast}\!h(y_{i})\mu_{x}(A_{i})\\
&  =\sum_{i}{}^{\ast}\!h(y_{i})\mu_{x_{0}}(A_{i})\dfrac{\mu_{x}(A_{i})}%
{\mu_{x_{0}}(A_{i})}.
\end{align*}
The hyperfinite set of weights $\left\{  ^{\ast}h(y_{i})\mu_{x_{0}}%
(A_{i})\right\}  $ forms an internal measure $\varphi_{h}$ on the indexed set
$\left\{  \dfrac{\mu_{\cdot}(A_{i})}{\mu_{x_{0}}(A_{i})}\right\}  $ of
internal harmonic functions on $U$. That is, for each $i$, the function
\[
z\mapsto\dfrac{\mu_{z}(A_{i})}{\mu_{x_{0}}(A_{i})}\text{, ~~~}z\in U
\]
is given the weight $^{\ast}h(y_{i})\mu_{x_{0}}(A_{i})$. Using the
construction in \cite{Loeb}, the measure $\varphi_{h}$ is converted to a
standard measure $L\left(  \varphi_{h}\right)  $, still on the same set of
internally harmonic functions. With an early but specific use of the
measurability of the standard part map applied to the mapping%
\[
\dfrac{\mu_{\cdot}(A_{i})}{\mu_{x_{0}}(A_{i})}\mapsto\operatorname{uccst}%
\left(  \dfrac{\mu_{\cdot}(A_{i})}{\mu_{x_{0}}(A_{i})}\right)  ,
\]
the now standard measure $L\left(  \varphi_{h}\right)  $ is moved to a
probability measure $P_{h}$ on $\mathcal{H}^{1}$.

For each $h\in\mathcal{H}^{1}$, for each $x\in W$,%
\[
h(x)=\int_{g\in\mathcal{H}^{1}}g(x)dP_{h}(g).
\]
That is, $P_{h}$ is a representing measure for $h$. If $h$ equals an affine
combination $\sum_{j}\alpha_{j}h_{j}$ of functions in $\mathcal{H}^{1}$, then
$P_{h}=\sum_{j}\alpha_{j}P_{h_{j}}$. A consequence of this fact, communicated
to the author by B. Fuchssteiner, is that the Fuchssteiner corollary in
\cite{Fuchssteiner} of a theorem of Cartier, Fell, and Meyer (see
\cite{Loeb-Boundary} ), shows that $P_{h}$ is the unique representing measure
for $h$ on the extreme points of $\mathcal{H}^{1}$.

\bigskip

The problem remains to connect this construction of representing measures with
the boundary $\partial W$ formed using our equivalence relation. For each
$h\in\mathcal{H}^{1}$, the internal measure we have constructed can also be
formed on $C$ using the weights $^{\ast}h(y_{i})\mu_{x_{0}}(A_{i}%
)\delta_{y_{i}}$. Using \cite{Loeb}, that measure can be transformed into a
standard probability measure on $C$, and then moved to $\partial W$ using the
standard part map. It is more important, however, to map $\partial W$ into
$\mathcal{H}^{1}$.

For each $z\in\partial W$, let $C_{z}$ denote the set of remote points in the
equivalence class forming $z$ but also in $C$. Every connected standard
neighborhood of $z$ contains points of $W$, so the nonstandard extension
contains points of $C$. It follows that $C_{z}$ is not empty. If $y\in C_{z}$
is in $A_{i}$ for some $i$, then it is associated with the harmonic function%
\[
\operatorname{uccst}\left(  \dfrac{\mu_{\cdot}(A_{i})}{\mu_{x_{0}}(A_{i}%
)}\right)  \in\mathcal{H}^{1}.
\]
If $y\in C_{z}$ is not in any $A_{i}$, we associate no function with $y$.

We now \textbf{assume} that for each $z\in\partial W$, we associate the same
element of $\mathcal{H}^{1}$ with each $y\in C_{z}$ for which we have
associated a function. This forms a map from all or part of $\partial W$ into
$\mathcal{H}^{1}$, and the representing measure on $\mathcal{H}^{1}$ can then
be viewed as a measure on $\partial W$.

\section{}

\end{document}